\documentclass{amsart}
\usepackage{amssymb}
\usepackage{amsmath}
\usepackage{a4}
\newtheorem{theorem}{Theorem}[section]
\newtheorem{lemma}[theorem]{Lemma}
\newtheorem{remark}[theorem]{Remark}
\newtheorem{definition}[theorem]{Definition}

\def\proof{\medbreak\noindent{\it Proof.\ }}
\def\endproof{\hfill\qedbox\smallbreak\noindent}
\def\qedbox{\hbox{$\rlap{$\sqcap$}\sqcup$}}
\makeatletter
  
  \@addtoreset{equation}{section}
 \makeatother

\begin{document}
\title[Completeness and the Osserman condition
  in signature $(2,2)$]
{Completeness, Ricci blowup, the Osserman and the conformal
 Osserman condition for Walker signature $(2,2)$ manifolds}%
\author[Brozos-V\'azquez et. al.]
{M. Brozos-V\'azquez$^1$,
 P. Gilkey$^2$, E. Garc\'{\i}a--R\'{\i}o${}^1$,\\
 and R. V\'azquez-Lorenzo${}^1$.}%
\begin{address}{B-V, G-R, V-L: Department of Geometry and Topology, Faculty of Mathematics,
University of Santiago de Compostela, Santiago de Compostela 15782, Spain.}\end{address}
\begin{email}{mbrozos@usc.es, xtedugr@usc.es, ravazlor@usc.es}\end{email}
\begin{address}{G: Mathematics Department, University of Oregon, Eugene Or 97403, USA}
\end{address}\begin{email}{gilkey@uoregon.edu}\end{email}

\begin{abstract}
Walker manifolds of signature $(2,2)$ have been used by many
authors to provide examples of Osserman and of  conformal Osserman
manifolds of signature $(2,2)$. We study questions of geodesic completeness and
Ricci blowup in this context.\newline
\end{abstract}

\maketitle

\section{Introduction}\label{sect-1}
We adopt the following notational conventions. Let
$\mathcal{M}:=(M,g)$ be a {\it pseudo-Riemannian manifold} of
signature $(p,q)$ and dimension $m=p+q$, let $\nabla$ be the {\it Levi-Civita
connection} of $\mathcal{M}$, let
$\mathcal{R}(x,y):=\nabla_x\nabla_y-\nabla_y\nabla_x-\nabla_{[x,y]}$
be the {\it curvature operator}, let
$\mathcal{J}(x):y\rightarrow\mathcal{R}(y,x)x$ be the {\it Jacobi
operator}, and let $S^\pm(\mathcal{M})$ be the {\it pseudo-sphere
bundles} of unit spacelike ($+$) and unit timelike ($-$) tangent
vectors.

Relating properties of the spectrum of the Jacobi operator $\mathcal{J}$ to the underlying geometry
of the manifold is an important area of investigation in recent years; we refer to \cite{GKVm02} for
a fuller discussion than is possible here. One says that
$\mathcal{M}$ is {\it pointwise spacelike} (resp. {\it timelike}) {\it Osserman} if the spectrum
of $\mathcal{J}$ is constant on
$S^+_P(\mathcal{M})$ (resp. on $S^-_P(\mathcal{M}))$ for every point $P\in M$. If $p>0$ and
$q>0$, these are equivalent concepts \cite{GKVV99} so we shall simply speak of a {\it pointwise
Osserman} manifold; in fact if $p>0$ and if $q>0$, one need only assume the spectrum of
$\mathcal{J}$ is bounded on $S^+(\mathcal{M})$ or on
$S^-(\mathcal{M})$
to ensure $\mathcal{M}$ is pointwise Osserman \cite{B06}. One replaces the word `pointwise' by
`globally' if the spectrum does not in fact depend on $P$.

The field began with a question raised by Osserman \cite{Os90} in
the Riemannian setting. Let $\mathcal{M}$ be a Riemannian
$2$-point homogeneous space, i.e. $\mathcal{M}$ is either flat or
is locally isometric to a rank $1$-symmetric space. Osserman noted
that the local isometries of $\mathcal{M}$ act transitively on the
bundle $S^+(\mathcal{M})$ and hence, in the notation adopted by
subsequent authors, $\mathcal{M}$ is globally Osserman. He
wondered if the converse held: is any globally Osserman Riemannian
manifold a local $2$-point homogeneous space? This question has
been called the {\it Osserman conjecture} and has been answered in
the affirmative by the work of Chi \cite{C88} and of Nikolayevsky
\cite{N04} except (possibly) in dimension 16 where the question is
still open. There is a similar classification result in the
Lorentzian setting. It is known \cite{BBG97,GKV97} that any
locally Osserman Lorentzian manifold has constant sectional
curvature. In the higher signature setting, such classification
results fail. There are, for example, Osserman manifolds that are
not  even locally affine homogeneous \cite{DGN05,GN05}.

Let $\{e_i\}$ be a local frame for the tangent bundle and let
$m=\dim (M)$. We set $g_{ij}:=g(e_i,e_j)$ and let $g^{ij}$ be the
inverse matrix. The {\it Ricci operator} $\rho$, the associated
{\it Ricci tensor} $\rho(\cdot,\cdot)$, the {\it scalar curvature}
$\tau$, the {\it Weyl conformal curvature} operator $\mathcal{W}$,
and the {\it conformal Jacobi operator} $\mathcal{J}_W$ are given
by:
\begin{eqnarray*}
&&\rho x:=\textstyle\sum_{ij}g^{ij}\mathcal{R}(x,e_i)e_j,\quad\rho(x,y):=g(\rho x,y),\quad
  \tau:=\sum_{ij}g^{ij}\rho(e_i,e_j),\\
&&\mathcal{W}(x,y)z=
\mathcal{R}(x,y)z+{\textstyle\frac1{(m-1)(m-2)}}\tau\{g(y,z)x-g(x,z)y\}\\
&&\qquad-{\textstyle\frac1{m-2}}\{g(x,z)\rho y+g(\rho x,z)y-g(y,z)\rho x-g(\rho y,z)x\},\\
&&\mathcal{J}_W(x):y\rightarrow\mathcal{W}(y,x)x\,.
\end{eqnarray*}
One says that $\mathcal{M}$ is {\it conformal Osserman} if for any
$P$ in $M$ one has that $\mathcal{J}_W$ has constant spectrum on
$S^+(\mathcal{M},P)$ for $p>0$ or, equivalently, on
$S^-(\mathcal{M},P)$ for $q>0$. This is a conformal notion
\cite{BBGS05}; if $\tilde g=e^\alpha g$ is conformally equivalent
to $g$, then $(M,g)$ is conformal Osserman if and only if
$(M,\tilde g)$ is conformal Osserman. Since any pointwise Osserman
manifold is necessarily Einstein, any pointwise Osserman manifold
is necessarily conformal Osserman.

We say that $\mathcal{M}$ is {\it geodesically complete} if all geodesics exist for all time. We say that
$\mathcal{M}$ exhibits {\it Ricci blowup} if there exists a geodesic $\gamma$ defined for $t\in[0,T)$ with
$T<\infty$ and if $\lim_{t\rightarrow T}|\rho(\dot\gamma,\dot\gamma)|=\infty$. Clearly if $\mathcal{M}$ exhibits
Ricci blowup, then it is geodesically incomplete and it can not be isometrically embedded in a geodesically complete
manifold.

In this paper, we will concentrate on signature $(2,2)$ where a great deal is known; the classification of Osserman
algebraic curvature tensors of signature $(2,2)$ is complete, see for example the discussion in \cite{BBZa01a}. Our
focus will be to relate the Osserman condition and the conformal Osserman condition, which are purely algebraic
conditions, to the global geometry of the manifold by studying questions of geodesic completeness and of Ricci
blowup. Here is a brief outline to the paper. In Section
\ref{sect-2}, we introduce the family of Walker manifolds we shall be considering and give their geodesic equations.
In Section
\ref{sect-3}, we show any strict Walker manifold is nilpotent Osserman and geodesically complete. In Section
\ref{sect-4}, we prove a result from the theory of ODEs which we will use subsequently to establish geodesic
incompleteness and Ricci blowup.

In the indefinite setting, the spectrum of a self-adjoint operator
does not determine the Jordan normal form of the operator. Let
$k\in\mathbb{R}$. There is a family of Walker manifolds of
signature $(2,2)$ whose Jacobi operator has eigenvalues
$\{0,4k,k,k\}$ but whose Jacobi operator is not diagonalizable
\cite{DGV06}. In Section \ref{sect-5} we show these manifolds
exhibit {\it Ricci blowup}. In Section \ref{sect-6}, we discuss
examples of conformal Osserman manifolds which exhibit various
eigenvalue structures following the discussion in \cite{BGVa05}.
We show that some of these manifolds are geodesically complete and
others exhibit Ricci blowup. Throughout this paper, we shall focus
on the global geometry of these manifolds and refer to previous
results in the literature for the corresponding algebraic features
of the curvature tensor.

\section{Signature $(2,2)$ Walker manifolds}\label{sect-2}

One says that a pseudo-Riemannian manifold $\mathcal{M}$ of signature $(2,2)$ is a {\it Walker
manifold} if it admits a parallel totally isotropic $2$-plane field. We refer to \cite{CGM05} for further
details. Such a manifold is locally isometric to an example of the following form:

\begin{definition}\label{defn-2.1}
\rm Let $(x_1,x_2,x_3,x_4)$ be coordinates on $\mathbb{R}^4$. Let
$\psi_{ij}(\vec x)=\psi_{ji}(\vec x)$ be smooth functions for $i,j=3,4$. Let
$\partial_i:=\partial_{x_i}$. Let
$\mathcal{M}:=(\mathbb{R}^4,g)$ where:
\begin{eqnarray*}
&&g(\partial_1,\partial_3)=g(\partial_2,\partial_4)=1,\quad
g(\partial_i,\partial_j)=\psi_{ij}\text{ for }i,j=3,4\,.
\end{eqnarray*}
\end{definition}

\begin{lemma}\label{lem-2.2}
Let $\psi_{ij/k}:=\partial_k\psi_{ij}$ and let $\mathcal{M}$ be as in Definition \ref{defn-2.1}.
Then the geodesic equations for $\mathcal{M}$ are given by:
\medbreak\quad$0=\ddot x_1+\dot x_1\dot x_3\psi_{33/1}+\dot x_1\dot x_4\psi_{34/1}+\dot x_2\dot
x_3\psi_{33/2}
       +\dot x_2\dot x_4\psi_{34/2}$\par\qquad$
       +\frac12\dot x_3\dot x_3(\psi_{33/3}+\psi_{34}\psi_{33/2}+\psi_{33}\psi_{33/1})
       $\par\qquad$
       +\dot x_3\dot x_4(\psi_{33/4}+\psi_{34}\psi_{34/2}+\psi_{33}\psi_{34/1})
       $\par\qquad$
       +\frac12\dot x_4\dot x_4(2\psi_{34/4}-\psi_{44/3}+\psi_{34}\psi_{44/2}+\psi_{33}\psi_{44/1})$,
       \smallbreak\quad
     $0=\ddot x_2+\dot x_1\dot x_3\psi_{34/1}+\dot x_1\dot x_4\psi_{44/1}+\dot x_2\dot x_3\psi_{34/2}
       +\dot x_2\dot x_4\psi_{44/2}$\par\qquad$
       +\frac12\dot x_3\dot x_3(2\psi_{34/3}-\psi_{33/4}+\psi_{44}\psi_{33/2}+\psi_{34}\psi_{33/1})
        $\par\qquad$
       +\dot x_3\dot x_4(\psi_{44/3}+\psi_{44}\psi_{34/2}+\psi_{34}\psi_{34/1})
       $\par\qquad$
       +\frac12\dot x_4\dot x_4(\psi_{44/4}+\psi_{44}\psi_{44/2}+\psi_{34}\psi_{44/1})$,
       \smallbreak\quad $0=\ddot x_3-\frac12\dot x_3\dot
x_3\psi_{33/1}-\dot x_3\dot x_4\psi_{34/1}
        -\frac12\dot x_4\dot x_4\psi_{44/1}$,
\smallbreak\quad $0=\ddot x_4-\frac12\dot x_3\dot x_3\psi_{33/2}-\dot x_3\dot x_4\psi_{34/2}
         -\frac12\dot x_4\dot x_4\psi_{44/2}$.
\end{lemma}
\proof We use the discussion in
\cite{GRVa05} to determine the Christoffel symbols; the Lemma then
follows: \medbreak\qquad$ \nabla_{\partial_1}\partial_3=
  \frac12\psi_{33/1}\partial_1+\frac12\psi_{34/1}\partial_2$,\quad$
\nabla_{\partial_1}\partial_4
       =\frac12\psi_{34/1}\partial_1+\frac12\psi_{44/1}\partial_2,
  $\par\qquad$
\nabla_{\partial_2}\partial_3
       =\frac12\psi_{33/2}\partial_1+\frac12\psi_{34/2}\partial_2$,\quad$
\nabla_{\partial_2}\partial_4
       =\frac12\psi_{34/2}\partial_1+\frac12\psi_{44/2}\partial_2,
$\par\qquad$ \nabla_{\partial_3}\partial_3=
\frac12(2\psi_{34/3}-\psi_{33/4}+\psi_{44}\psi_{33/2}+\psi_{34}\psi_{33/1})\partial_2
$\par\qquad\qquad\quad\,$
       +\frac12(\psi_{33/3}+\psi_{34}\psi_{33/2}+\psi_{33}\psi_{33/1})\partial_1
-\frac12\psi_{33/1}\partial_3-\frac12\psi_{33/2}\partial_4,
$\par\qquad$ \nabla_{\partial_3}\partial_4
       =\frac12(\psi_{33/4}+\psi_{34}\psi_{34/2}+\psi_{33}\psi_{34/1})\partial_1
$\par\qquad\qquad\quad\,$
+\frac12(\psi_{44/3}+\psi_{44}\psi_{34/2}+\psi_{34}\psi_{34/1})\partial_2
-\frac12\psi_{34/1}\partial_3-\frac12\psi_{34/2}\partial_4,
$\par\qquad$ \nabla_{\partial_4}\partial_4
      =\frac12(2\psi_{34/4}-\psi_{44/3}+\psi_{34}\psi_{44/2}+\psi_{33}\psi_{44/1})\partial_1
$\par\qquad\qquad\quad\,$
       +\frac12(\psi_{44/4}+\psi_{44}\psi_{44/2}+\psi_{34}\psi_{44/1})\partial_2
       -\frac12\psi_{44/1}\partial_3-\frac12\psi_{44/2}\partial_4$.\hfill\endproof

\medbreak We specialize this to the following situation. Assertion (1) in the following Lemma is a consequence of
Lemma
\ref{lem-2.2} while Assertion (2) in the following Lemma follows by specializing results of \cite{DGV06a}:
\begin{lemma}\label{lem-2.3}
Let $\mathcal{M}:=(\mathbb{R}^4,g)$ where the non-zero components
of $g$ are given by
$g(\partial_1,\partial_3)=g(\partial_2,\partial_4)=1$ and
$g(\partial_3,\partial_4)=\psi_{34}(x_1,x_2,x_3,x_4)$.
Then:\begin{enumerate} \item The geodesic equations are:
\par\qquad $0=\ddot x_1+\dot x_1\dot x_4\psi_{34/1}+\dot x_2\dot x_4\psi_{34/2}+\dot x_3\dot
x_4\psi_{34}\psi_{34/2}+\dot x_4\dot x_4\psi_{34/4}$,\par\qquad
$0=\ddot x_2+\dot x_1\dot x_3\psi_{34/1}+\dot x_2\dot x_3\psi_{34/2}+\dot x_3\dot x_3\psi_{34/3}
+\dot x_3\dot x_4\psi_{34}\psi_{34/1}$,\par\qquad
$0=\ddot x_3-\dot x_3\dot x_4\psi_{34/1}$,\qquad
$0=\ddot x_4-\dot x_3\dot x_4\psi_{34/2}$.
\item The non-zero components of the Ricci tensor are:
\par\qquad
$\rho_{13}=\rho_{24}=\frac12\psi_{34/12}$,
\qquad$\rho_{14}=\frac12\psi_{34/11}$,
\qquad$\rho_{23}=\frac12\psi_{34/22}$,\par\qquad
$\rho_{33}=\frac12\{-\psi_{34/2}^2+2\psi_{34/23}\}$,
\qquad$\rho_{44}=\frac12\{-\psi_{34/1}^2+2\psi_{34/14}\}$,\par\qquad
$\rho_{34}=\textstyle\frac12\{\psi_{34/1}\psi_{34/2}
      +2\psi_{34}\psi_{34/12}-\psi_{34/13}-\psi_{34/24}\}$.
\end{enumerate}\end{lemma}

\section{Strict Walker manifolds}\label{sect-3}

One says that $\mathcal{M}$ is {\it nilpotent Osserman} if $\mathcal{J}(x)$ has only the
eigenvalue $0$ or, equivalently, if $\mathcal{J}(x)^{\dim(M)}=0$ for any tangent vector $x$.
One says that $\mathcal{M}$ is a {\it strict Walker manifold} if $\psi_{ij}=\psi_{ij}(x_3,x_4)$.
Our first result deals with such examples:
\begin{theorem}\label{thm-3.1}
Let $\mathcal{M}$ be as in Definition \ref{defn-2.1} with $\psi_{ij}=\psi_{ij}(x_3,x_4)$.
Then $\mathcal{M}$ is nilpotent Osserman and geodesically
complete.
\end{theorem}

\proof By Lemma \ref{lem-2.2}, the geodesic equations in this setting are:
\begin{eqnarray*}
&&0=\textstyle\ddot
x_1+\frac12\dot x_3\dot x_3\psi_{33/3}
       +\frac12\dot x_4\dot x_4(2\psi_{34/4}-\psi_{44/3})
       +\dot x_3\dot x_4\psi_{33/4},\\
&&0=\textstyle\ddot x_2+\frac12\dot x_3\dot
x_3(2\psi_{34/3}-\psi_{33/4})
       +\frac12\dot x_4\dot x_4\psi_{44/4}
       +\dot x_3\dot x_4\psi_{44/3},\\
&&0=\textstyle\ddot x_3,\quad\quad 0=\ddot x_4\,.
\end{eqnarray*}
The final two equations can be solved to yield $x_3=a+bt$ and $x_4=c+dt$. The first two equations
then have the form $\ddot x_1=f_1(t)$ and $\ddot x_2=f_2(t)$ which can be solved. Thus geodesics
extend for infinite time and $\mathcal{M}$ is geodesically complete. One computes easily that
$\mathcal{J}(x):\operatorname{Span}\{\partial_3,\partial_4\}\rightarrow
\operatorname{Span}\{\partial_1,\partial_2\}\rightarrow0$
for any tangent vector $x$ and hence $\mathcal{J}(x)$ is nilpotent. This shows $\mathcal{M}$ is
nilpotent Osserman; we refer to the discussion in \cite{GIZ03} for further details.\endproof

\section{A result from ODEs}\label{sect-4}
Before continuing our investigations further, we shall need the following result:
\begin{lemma}\label{lem-4.1}
Let $f(t)$ satisfy $\ddot f(t)=\Xi(\dot f,f)$ with $f(0)=1$ and
$\dot f(0)=1$ and maximal domain $[0,T)$. Assume
$\Xi(x,y)\ge\varepsilon x^ay^b$ for $x\ge1$ and $y\ge1$ where $2a+b\ge3$ and $\varepsilon>0$.
Then $T<\infty$, $\lim_{t\rightarrow T}\ddot f(t)=\infty$, and
$\lim_{t\rightarrow T}\dot f(t)=\infty$.
\end{lemma}

\proof Since $\ddot f$ is positive, $f$ and $\dot f$ are monotonically increasing. Suppose $T<\infty$ but
$\lim_{t\rightarrow T}\dot f<\infty$. Then
$\dot f$  is bounded and hence $f$ is bounded as well. Thus
$\lim_{t\rightarrow T}\dot f=\dot f_T$ and
$\lim_{t\rightarrow T}f=f_T$ exist and are finite and the fundamental
theorem of ODEs shows
$[0,T)$ is not the maximal domain of the function $f$. Thus if $T$ is finite, $\lim_{t\rightarrow T}\dot
f(t)=\infty$. Hence $\lim_{t\rightarrow T}\ddot f(t)=\infty$ as
well and the Lemma holds.

To complete the proof, we suppose that $T=\infty$ and argue for a contradiction. Without loss of
generality, we assume $\varepsilon<1$. Let
$t_1:=0$ and let
$t_{n+1}:=t_n+\frac3{\varepsilon n^2}$. We wish to show
$f(t_n)\ge n$ and $\dot f(t_n)\ge n^2$. As this holds for $n=1$, we proceed by induction on $n$.
Because $\dot f$ is
monotonically increasing,
$$f(t_{n+1})\ge f(t_n)+\dot f(t_n)\textstyle\frac3{\varepsilon n^2}
\ge n+\frac{3n^2}{\varepsilon n^2}\ge n+1\,.$$
To prove the second estimate, we use the mean value theorem to choose $s$ with
$s\in[t_n,t_{n+1}]$ so
$\dot f(t_{n+1})=\dot f(t_n)+\textstyle\frac3{\varepsilon n^2}\ddot f(s)$. We may estimate that:
\begin{eqnarray*}
\ddot f(s)&=&\Xi(\dot f(s),f(s))\ge\varepsilon\dot f(s)^af(s)^b
\ge\varepsilon \dot f(t_n)^af(t_n)^b
\ge \varepsilon n^{2a+b}\ge\varepsilon n^3\,.
\end{eqnarray*}
Consequently
$\dot f(t_{n+1})\ge n^2+\textstyle\frac3{\varepsilon n^2}\varepsilon n^3\ge n^2+3n\ge(n+1)^2$.
The desired contradiction follows as $\lim_{n\rightarrow\infty}t_n<\infty$.
\endproof

\section{Non-diagonalizable Jacobi operators}\label{sect-5}

In signature $(2,2)$, the eigenvalue structure of the Jacobi operator does not determine the
operator up to conjugacy; one must instead consider the Jordan normal form. Theorem \ref{thm-3.1}
shows any strict Walker manifold of signature $(2,2)$ is nilpotent Osserman. However there are Walker manifolds of
signature $(2,2)$ which are Osserman but not nilpotent and whose Jacobi operators are not
diagonalizable.

\begin{theorem}\label{thm-5.1}
Let $\mathcal{M}$ be given by Definition \ref{defn-2.1} with
$\psi_{33}=4kx_1^2-\textstyle\frac1{4k}f(x_4)^2$,
$\psi_{44}=4kx_2^2$, and $\psi_{34}=4kx_1x_2+x_2f(x_4)-\frac1{4k}\dot
f(x_4)$ where $f=f(x_4)$ is non-constant and $k\neq 0$. Then
$\mathcal{M}$ is Osserman with eigenvalues $\{0,4k,k,k\}$ and the
Jacobi operators are diagonalizable at $P$ if and only if
\smallbreak\centerline{$24kf(x_4)\dot f(x_4)x_2-12k\ddot f(x_4)x_1
+3f(x_4)\ddot f(x_4)+4\dot f(x_4)^2=0$.}\smallbreak\noindent
Furthermore $\mathcal{M}$ is geodesically incomplete and can not
be embedded isometrically in a geodesically complete manifold.
\end{theorem}

\proof We refer to \cite{DGV06} for the proof that $\mathcal{M}$
is Osserman with the indicated Jordan normal form. We must show
$\mathcal{M}$ is geodesically incomplete. Since $f$ is
non-constant, we may choose $\xi_4$ so $f(\xi_4)\ne0$ and $\dot
f(\xi_4)\ne0$. Choose $\xi_1$ so $16k^2\xi_1^2=f(\xi_4)^2$;
normalize the choice of sign so $k\xi_1>0$. As an ansatz, we set
$x_1=\xi_1$, $x_2=0$, and $x_4=\xi_4$ to be constant. This implies
$\psi_{33}=0$. The geodesic equations in $\ddot x_1$, $\ddot x_2$,
and $\ddot x_4$ given by Lemma \ref{lem-2.2} then become $\ddot
x_1=\ddot x_2=\ddot x_4=0$ which are satisfied.  The remaining
geodesic equation is $0=\ddot x_3-4k\xi_1\dot x_3\dot x_3$. We can solve this equation by setting
$$x_3=-\textstyle\frac1{4k\xi_1}\ln(1-t),\
  \dot x_3=\frac1{4k\xi_1}(1-t)^{-1},\
  \ddot x_3=\frac1{4k\xi_1}(1-t)^{-2}=4k\xi_1\dot x_3\dot x_3\,.$$
This is defined for $t\in(-\infty,1)$ and we have $\lim_{t\rightarrow1}4k\xi_1x_3=\infty$. In
particular,
$\mathcal{M}$ is geodesically incomplete.

Since
$\mathcal{M}$ is Einstein, it does not exhibit Ricci blowup. Instead we use a different argument
to show $\mathcal{M}$ is essentially incomplete.
Let $\{e_1,e_2,e_3,e_4\}$ be a parallel frame along $\gamma$ with
$e_i(0)=\partial_i$.  The argument used to establish Lemma \ref{lem-2.2} shows:
$$\begin{array}{ll}
   \nabla_{\partial_3}\partial_1=4k\xi_1\partial_1,&
\nabla_{\partial_3}\partial_2=4k\xi_1\partial_2,\\
   \nabla_{\partial_3}\partial_3=-4k\xi_1\partial_3,&
  \nabla_{\partial_3}\partial_4=-2\xi_1\dot f(\xi_4)\partial_1-4k\xi_1\partial_4\,.
\end{array}$$
 Consequently
\begin{eqnarray*}
&&e_1(x_3)=e^{-4k\xi_1x_3}\partial_1,\quad e_2(x_3)=e^{-4k\xi_1x_3}\partial_2,\quad
e_3(x_3)=e^{4k\xi_1x_3}\partial_3,\\
&&e_4(x_3)=\textstyle\frac1{4k}\dot
f(\xi_4)(e^{4k\xi_1x_3}-e^{-4k\xi_1x_3})\partial_1+e^{4k\xi_1x_3}\partial_4\,.
\end{eqnarray*}
Since $R(\partial_1,\partial_3,\partial_3,\partial_4)=0$, since
$R(\partial_1,\partial_3,\partial_3,\partial_1)=4k$, since $\dot f(\xi_4)\ne0$, and since
$4k\xi_1x_3(t)\rightarrow\infty$ as
$t\rightarrow1$,
$\mathcal{M}$ is seen to be essentially incomplete as:
\begin{eqnarray*}
&&\lim_{t\rightarrow 1}R(e_1,e_3,e_3,e_4)
=\lim_{t\rightarrow 1}\{\textstyle\frac1{4k}\dot f(\xi_4)
(e^{4k\xi_1x_3}-e^{-4k\xi_1x_3})\}e^{4k\xi_1x_3}4k\\
&=&\lim_{t\rightarrow 1}\dot f(\xi_4)(e^{8k\xi_1x_3}-1)=\pm\infty\,.\qquad\qedbox
\end{eqnarray*}

\section{Conformal Osserman manifolds}\label{sect-6}

Let $\operatorname{Spec}_W$ denote the spectrum of the conformal
Jacobi operator and let $m_\lambda$ denote the minimal polynomial
of the conformal Jacobi operator for a conformal Osserman
manifold. We refer to \cite{BGVa05} for the proof of:

\begin{theorem}\label{thm-6.1}
Let $\mathcal{M}$ be as in Definition \ref{defn-2.1} where
$\psi_{33}=\psi_{44}=0$. With the following choices of
$\psi_{34}$, $\mathcal{M}$ is conformal Osserman and:
\begin{enumerate}
\item The Jordan normal form does not change from point to point:
\begin{enumerate}
\item\hglue -.1truecm $\psi_{34}=x_1^2-x_2^2\Rightarrow
m_\lambda=\lambda(\lambda^2-\frac14)$ and
$\operatorname{Spec}_W=\{0,0,\pm\frac12\}$. \item\hglue -.1truecm
$\psi_{34}=x_1^2+x_2^2\Rightarrow
m_\lambda=\lambda(\lambda^2+\frac14)$ and
$\operatorname{Spec}_W=\{0,0,\pm\frac{\sqrt{-1}}2\}$. \item\hglue
-.1truecm $\psi_{34}=x_1x_4+x_3x_4\Rightarrow m_\lambda=\lambda^2$
and $\operatorname{Spec}_W=\{0\}$. \item\hglue -.1truecm
$\psi_{34}=x_1^2\Rightarrow m_\lambda=\lambda^3$  and
$\operatorname{Spec}_W=\{0\}$.
\end{enumerate}
\item\hglue -.1truecm  $\operatorname{Spec}_W=\{0\}$ but Jordan normal form changes from point to
point.
\begin{enumerate}
\item\hglue -.1truecm $\psi_{34}=x_2x_4^2+x_3^2x_4\Rightarrow
m_\lambda =\left\{\begin{array}{lll}
   \lambda^3&\text{if}&x_4\ne0,\\
   \lambda^2&\text{if}&x_4=0,\ x_3\ne0,\\
   \lambda&\text{if}&x_3=x_4=0\,.\end{array}\right.$
\item\hglue -.1truecm $\psi_{34}=x_2x_4^2+x_3x_4\Rightarrow
m_\lambda =\left\{\begin{array}{lll}
\lambda^3&\text{if}&x_4\ne0,\\
\lambda^2&\text{if}&x_4=0\,.\end{array}\right.$
\item\hglue -.1truecm$\psi_{34}=x_1x_3^2\Rightarrow
m_\lambda =\left\{\begin{array}{lll}
\lambda^3&\text{if}&x_3\ne0,\\
\lambda&\text{if}&x_3=0.\end{array}\right.$
\item\hglue -.1truecm$\psi_{34}=x_1x_3+x_2x_4
\Rightarrow m_\lambda =\left\{\begin{array}{lll}
\lambda^2&\text{if}&x_1x_3+x_2x_4\ne0,\\
\lambda&\text{if}&x_1x_3+x_2x_4=0\,.
\end{array}\right.$
\end{enumerate}
\item\hglue -.1truecm The eigenvalues can change from point to point:
\begin{enumerate}
\item\hglue -.1truecm
$\psi_{34}=x_1^4+x_1^2-x_2^4-x_2^2\Rightarrow
\operatorname{Spec}_W=\{0,0,\pm{\textstyle\frac12}\sqrt{(6x_1^2+1)(6x_2^2+1)}\}$.
\item\hglue -.1truecm
$\psi_{34}=x_1^4+x_1^2+x_2^4+x_2^2\Rightarrow
\operatorname{Spec}_W=\{0,0,\pm{\textstyle\frac12}\sqrt{-(6x_1^2+1)(6x_2^2+1)}\}$.
\item\hglue -.1truecm
$\psi_{34}=x_1^3-x_2^3\Rightarrow\operatorname{Spec}_W=\{0,0,\pm{\textstyle\frac32}\sqrt{x_1x_2}\}
$.
\end{enumerate}\end{enumerate}
\end{theorem}

\begin{remark}\label{rmk-6.2}\rm In revisiting the manifolds of Theorem \ref{thm-6.1} whilst
writing this paper, we made some geometrical observations that, although not directly in the
focus of this paper, never the less illustrate why they form a rich geometrical family that
it is important to study. Recall that a manifold is said to be {\it curvature homogeneous} if
given $P,Q\in M$, there is an isometry
$\phi:T_PM\rightarrow T_QM$ so $\phi^*R_Q=R_P$ -- i.e.  the curvature
tensor ``looks the same at any point of the manifold''; we refer to
\cite{DGN05} for further details. Of the manifolds in Theorem \ref{thm-6.1}, only the
manifold of 1-d) is curvature homogeneous. Let
$R_\Lambda:\Lambda^2T^*M\rightarrow\Lambda^2T^*M$ be the curvature operator. If
$\mathcal{M}$ is the manifold of 1-a) or 1-b) above, then $\rho^2$ is a
multiple of the identity at a point if and only if $x_1=x_2=0$. If $\mathcal{M}$ is as in
1-c), then $\mathcal{M}$ is Ricci flat at a point if and only if $x_4^2=2$; for 2-a) and
2-b), this happens if and only if
$x_4=0$; for 2-c), this happens if and only if $x_3=0$. For 2-d),
$R_\Lambda^2=0$ if and only if $x_3x_4=1$. The eigenvalues of
$R_\Lambda$ change from point to point for the manifolds of Theorem
3-a), 3-b), and 3-c). These observations show that none of these manifolds are curvature
homogeneous. On the other hand, a rather more delicate argument shows Example 1-d) is
curvature homogeneous. Further details concerning these matters will be forthcoming in a
subsequent article.  We note finally that Derdzinski \cite{D00} showed a $4$-dimensional
Riemannian manifold is curvature homogeneous if and only if $R_\Lambda$ has constant
eigenvalues; furthermore if such a manifold is Einstein, then it is locally symmetric.
In Example 2-d), $\operatorname{Spec}\{R_\Lambda\}=\{0\}$ but the manifold is not curvature
homogeneous. Thus this result of Derdzinski fails in signature $(2,2)$; we refer to
Derdzinski \cite{D02} for additional results in this direction.
\end{remark}

We study the global geometry of the manifolds of Theorem \ref{thm-6.1}:
\begin{theorem}\label{thm-6.3}
Of the manifolds in Theorem \ref{thm-6.1}, only $\psi_{34}=x_1x_4+x_3x_4$ defines
 a geodesically complete manifold;
the remaining tensors $\psi$ define manifolds which exhibit Ricci blowup and which therefore can
not be embedded isometrically in a geodesically complete manifold.
\end{theorem}

\proof Suppose first that for $x_1\ge1$ one has:
$$
\psi_{34/1}=p(x_1)\ge x_1,\quad \psi_{34/11}\ge1,\quad \psi_{34/4}=0
\,.$$
This is the case for the warping functions of (1a), (1b), (1d), (3a), (3b), and (3c). We set
$x_2(t)=0$, $x_3(t)=0$ and
$x_4(t)=-t$. The geodesic equations given in Lemma \ref{lem-2.3} then become:
$$
\ddot x_1-\dot x_1p(x_1)=0,\quad
 \ddot x_2=0,\quad\ddot x_3=0,\qquad\ddot x_4=0\,.
$$
This yields a consistent set of equations with $\ddot x_1=\dot
x_1p(x_1)$. By Lemma \ref{lem-4.1}, $\lim_{t\rightarrow T}\dot
x_1(t)=\infty$ for $T$ finite. By Lemma \ref{lem-2.3},
\begin{eqnarray*}
\rho(\dot\gamma,\dot\gamma)&=&\dot x_1\dot x_1\rho_{11}+2\dot x_1\dot x_4\rho_{14}
        +\dot x_4\dot x_4\rho_{44}
=-\dot x_1\psi_{34/11}-\textstyle\frac12\psi_{34/1}^2\,.
\end{eqnarray*}
As $\psi_{34/11}\ge1$, $\lim_{t\rightarrow T}\rho(\dot\gamma,\dot\gamma)=-\infty$  and these manifolds exhibit Ricci
blowup.

 Let $\psi_{34}=x_1x_3+x_2x_4$ be as in (2d). The geodesic equations
are: \medbreak\quad$0=\ddot x_1+\dot x_1\dot x_4x_3+\dot x_2\dot
x_4x_4+\dot x_3\dot x_4(x_1x_3+x_2x_4)x_4+\dot x_4\dot x_4x_2$,
\smallbreak\quad$0=\ddot x_2+\dot x_1\dot x_3x_3+\dot x_2\dot
x_3x_4+\dot x_3\dot x_3x_1 +\dot x_3\dot x_4(x_1x_3+x_2x_4)x_3$,
\smallbreak\quad $0=\ddot x_3-\dot x_3\dot x_4x_3$,\qquad
 $0=\ddot x_4-\dot x_3\dot x_4x_4$.
\medbreak\noindent We start with initial conditions $x_3(0)=\dot x_3(1)=x_4(1)=\dot x_4(1)=1$. Symmetry implies that
$x_3(t)=x_4(t)=h(t)$ where $h$ satisfies
$\ddot h(t)=\dot h(t)\dot h(t)h(t)$.
Lemma \ref{lem-4.1} now shows $\dot h\rightarrow\infty$ at finite time so $\mathcal{M}$ is incomplete. Furthermore,
we use Lemma \ref{lem-2.3} to see that $\rho_{ij}=0$ for $i,j\ne3,4$ and thereby show
$\mathcal{M}$ exhibits Ricci blowup by computing:
\begin{eqnarray*}
\rho(\dot\gamma,\dot\gamma)&=&\dot x_3^2\rho_{33}+\dot x_4^2\rho_{44}+2\dot x_3\dot x_4\rho_{34}
=\dot x_3^2\{\rho_{33}+\rho_{44}+2\rho_{34}\}\\
&=&\textstyle\dot x_3^2\{-\frac12x_4^2-\frac12x_3^2+x_3x_4-2\}=-2\dot x_3^2\,.
\end{eqnarray*}

Let $\psi_{34}=x_1x_3^2$ be the warping function of (2c). The
final two geodesic equations become $\ddot x_3=\dot x_3\dot x_4
x_3^2$ and $\ddot x_4=0$. Setting $x_4=t$ then yields the equation
$\ddot x_3=\dot x_3x_3^2$ and thus by Lemma \ref{lem-4.1}, $\dot
x_3(t)\rightarrow\infty$ as $t\rightarrow T$ for $T<\infty$. All
the components of the Ricci tensor vanish except $\rho_{34}$ and
$\rho_{44}$. Since $x_3(t)\ge1$, we show $\mathcal{M}$ exhibits
Ricci blowup by computing:
\begin{eqnarray*}
\lim_{t\rightarrow T}\rho(\dot\gamma,\dot\gamma)&=&
\lim_{t\rightarrow T}\{-2\dot x_3\dot x_4x_3-{\textstyle\frac12}\dot x_4\dot x_4x_3^4\}
\le\lim_{t\rightarrow T}\{-2\dot x_3\}=-\infty\,.
\end{eqnarray*}

Suppose that $\psi_{34}=x_2x_4^2+x_3^2x_4$ or that
$\psi_{34}=x_2x_4^2+x_3x_4$ are the warping functions of Theorem
\ref{thm-6.1} (2a) and Theorem \ref{thm-6.1} (2b). The final two
geodesic equations become $0=\ddot x_3$ and $\ddot x_4=\dot
x_3\dot x_4x_4^2$. We take $x_3=t$ so $\ddot x_4=\dot x_4x_4^2$.
Let $x_4(0)=\dot x_4(0)=1$.
Thus $\lim_{t\rightarrow T}\dot x_4=\infty$ at some finite time and $x_4(t)\ge1$ for all $t$.
Only $\rho_{33}$ and $\rho_{34}$ are non-zero. Thus we may show
$\mathcal{M}$ exhibits Ricci blowup by computing:
\begin{eqnarray*}
\lim_{t\rightarrow T}\rho(\dot\gamma,\dot\gamma)&=&
\lim_{t\rightarrow T}\left\{-\textstyle\frac12\dot x_3\dot x_3x_4^4-2\dot x_3\dot x_4x_4\right\}
\le\lim_{t\rightarrow T}\{-2\dot x_4\}=-\infty\,.
\end{eqnarray*}

Finally let $\psi_{34}=x_1x_4+x_3x_4$ be the warping function of
(1c). The
geodesic equations in the last two variables are $\ddot x_3-\dot x_3\dot x_4x_4=0$ and
$\ddot x_4=0$ so
$$
x_4=a+bt\quad\text{and}\quad\dot x_3=ce^{b(at+b\frac12t^2)}\,.
$$
We integrate this equation to determine $x_3$. As the equation for
$\ddot x_2$ takes the form $\ddot x_2+F(x_1,x_3,x_4,\dot x_1,\dot
x_3,\dot x_4)=0$, it poses no difficulty and only task is to
determine $x_1$. The equation for $x_1$ takes the form:
$$\ddot x_1+\dot x_1\dot x_4x_4+\dot x_4\dot x_4x_1+\dot x_4\dot x_4 x_3=0\,.$$
By rescaling the geodesic parameter, we see that there are really
only two cases to be considered. These are $x_4=a$ and $x_4=t$.
If $x_4=a$, we get the equation $\ddot x_1=0$ which has linear
solutions. If $x_4=t$, we get the equation
$$\ddot x_1+t\dot x_1 +x_1=\alpha(t)$$
for suitably chosen $\alpha$. We set $x_1:=fe^{-\frac12t^2}$ to reduce the order of the equation:
\begin{eqnarray*}
&&\dot x_1=(\dot f-tf)e^{-\frac12t^2},\quad
\ddot x_1=(\ddot f-2t\dot f+t^2f-f)e^{-\frac12t^2},\\
&&\ddot x_1+t\dot x_1+x_1=(\ddot f-2t\dot f+t^2f-f+t\dot f-t^2f+f)e^{-\frac12t^2}\\
&&\qquad\quad=(\ddot f-t\dot f)e^{-\frac12t^2}=\alpha(t)\,.
\end{eqnarray*}
Setting $f_1:=\dot f$ then leads to an equation of the form
$\dot f_1-tf_1=\alpha_1(t)$ for suitably chosen $\alpha_1$. Setting $f_1=f_2e^{\frac12t^2}$ then
yields
$$\dot f_1-tf_1=(\dot f_2+tf_2-tf_2)e^{\frac12t^2}=\alpha_1(t)$$
which leads to the equation $\dot f_2=\alpha(t)$. This equation
can be solved for all time; it now follows that $\mathcal{M}$ is geodesically complete.\endproof

\section*{Acknowledgments} Research of M. Brozos-V\'azquez and of P. Gilkey
was partially supported by the Max Planck Institute for Mathematics in
the Sciences (Germany). Research of M. Brozos-V\'azquez was also
partially supported by FPU grant and by project MTM2005-08757-C04-01
(Spain). Research of E. Garc\'{\i}a--R\'{\i}o and of R. V\'azquez-Lorenzo
was partially supported by BFM2003-02949 (Spain).

  \end{document}